\documentclass[11pt]{article}
\usepackage{amssymb,amsfonts,amsmath, psfrag, amscd, cite}
\usepackage{graphicx}
\newtheorem{theo}{Theorem}

\newtheorem{defi}[theo]{Definition}
\newtheorem{lemm}[theo]{Lemma}
\newtheorem{rema}[theo]{Remark}

\makeatletter \@addtoreset{equation}{section}
\@addtoreset{theo}{section} \makeatother

\def\qed{\hfill \rule{4pt}{7pt}}
\def\pf{\noindent {\it Proof.} }

\begin{document}

\begin{center}
{\bf\LARGE Monochromatic and Heterochromatic\\[6pt] Subgraph
Problems\\[6pt] in a Randomly Colored Graph\footnote{Supported by NSFC, PCSIRT and the ``973"
program.}}
\end{center}
\vskip 1cm
\begin{center}
Xueliang Li and Jie Zheng\\
\vskip 0.5cm

Center for Combinatorics, LPMC-TJKLC\\
Nankai University\\
Tianjin 300071, P.R. China

Email: x.li@eyou.com, jzheng@eyou.com
\end{center}

{\bf Abstract.} Let $K_n$ be the complete graph with $n$ vertices
and $c_1, c_2, \cdots, c_r$ be $r$ different colors. Suppose we
randomly and uniformly color the edges of $K_n$ in  $c_1, c_2,
\cdots, c_r$. Then we get a random graph, denoted by
$\mathcal{K}_n^r$. In the paper, we investigate the asymptotic
properties of several kinds of monochromatic and heterochromatic
subgraphs in $\mathcal{K}_n^r$. Accurate threshold functions in
some cases are also obtained.



\section{Introduction}
The study of random graphs was begun by P. Erd\"{o}s and A.
R\'{e}nyi in the 1960s \cite{erdo1, erdo2, erdo3} and now has a
comprehensive literature \cite{boll, jans}.

The most frequently encountered probabilistic model of random
graph is $\mathcal{G}_{n,p(n)}$, where $0\leq p(n)\leq 1$.  It
consists of all graphs with vertex set $V=\{1,2,\cdots,n\}$ in
which the edges are chosen independently and with probability
$p(n)$. As $p(n)$ goes from zero to one the random graph
$\mathcal{G}_{n,p(n)}$ evolves from empty to full.

P. Erd\"{o}s and A. R\'{e}nyi discovered that for many natural
properties A of graphs there was a narrow range in which
$Pr[\mathcal{G}_{n,p(n)}\mbox{ has property A}]$ moves from near
zero to near one. So we introduce the following important
definition (\cite{spen}, page 14).

\begin{defi}\label{thre}
A function $p(n)$ is a threshold function for property A if the
following two conditions are satisfied:
\begin{enumerate}
\item If $p'(n)\ll p(n)$, then $\lim_{n\to\infty}Pr[\mathcal{G}_{n,p'(n)}\mbox{ has property A}]=0$.\\
\item If $p'(n)\gg p(n)$, then
$\lim_{n\to\infty}Pr[\mathcal{G}_{n,p'(n)}\mbox{ has property
A}]=1$.
\end{enumerate}
\end{defi}

In general, if $Pr[\mathcal{G}_{n,p(n)}\mbox{ has property A}]\to
0$, we say {\it almost no} $\mathcal{G}_{n,p(n)}$ has property A.
Conversely, if $Pr[\mathcal{G}_{n,p(n)}\mbox{ has property A}]\to
1$, we say {\it almost every} $\mathcal{G}_{n,p(n)}$ has property
A.

In this article, we introduce the following probabilistic model of
random graphs. Let $K_n$ be the complete graph with vertex set
$V=\{1,2,\cdots,n\}$ and $c_1,c_2,\cdots,c_r$ be $r=r(n)$
different colors. We now send $c_1,c_2,\cdots,c_r$ to the edges of
$K_n$ randomly and equiprobably, which means each edge is colored
in $c_i(1\leq i\leq r)$ with probability $\frac{1}{r}$. Thus we
get a random graph $\mathcal{K}_n^r$. The probability space
$(\Omega,\mathcal{F},\mathcal{P})$ of $\mathcal{K}_n^r$ has a
simple form: $\Omega$ has $r^{n\choose 2}$ elements and each one
has probability $\frac{1}{r^{n\choose 2}}$ to appear.

The subgraph of $\mathcal{K}_n^r$ with vertices $1,2,\cdots,n$ and
the edges that have color $c_i$ is denote by $\mathcal{G}_i$.
Obviously, it is just the random graph $\mathcal{G}_{n,p(n)}$
(\cite{boll}, page 34), where $p(n)=\frac{1}{r}$.

Matching, clique and tree are three kinds of important subgraphs.
As to their definitions, please refer to \cite{bond}. A
$k$-matching is a matching of $k$ independent edges. A $k$-clique
is a clique of $k$ vertices. Similar, a $k$-tree is a tree of $k$
vertices. In a $k$-matching ($k$-clique, $k$-tree), if all of the
edges are in a same color, we call it a monochromatic $k$-matching
($k$-clique, $k$-tree); On the other hand, if any two of edges are
of different colors, we call it a heterochromatic $k$-matching
($k$-clique, $k$-tree).

Having a monochromatic $k$-matching, $k$-clique or $k$-tree or a
heterochromatic $k$-matching, $k$-clique or $k$-tree are all
properties of $\mathcal{K}_n^r$. We want to investigate these
properties and obtain the threshold functions for them. Two
properties will be especially demonstrated: monochromatic
$k$-matching and heterochromatic $k$-matching. For the others, the
methods are similar and we list the results in Section 4.

\section{Monochromatic $k$-Matchings in $\mathcal{K}_n^r$}
Let $k$ be an integer. Obviously, in $\mathcal{K}_n^r$, there are
altogether
\begin{equation*}
q=\frac{{n\choose 2}{{n-2}\choose 2}\cdots{{n-2k+2}\choose 2}}{k!}
\end{equation*}
sets of $k$ independent edges. Arrange them in an order and the
$i$-th one is denoted by $M_i$.

Let $A_i$ be the event that the edges in $M_i$ are monochromatic
and $X_i$ be the indicator variable for $A_i$. That is,
\begin{equation}
X_i= \left\{
\begin{array}{ll}
1 & \mbox{if } A_i \mbox{ happens},\\
0 & \mbox{otherwise}.
\end{array}
\right.
\end{equation}

Then the random variable
$$X=X_1+X_2+\cdots+X_q$$
denotes the number of monochromatic $k$-matchings in
$\mathcal{K}_n^r$.

For each $1\leq i\leq q$,
\begin{equation*}
E(X_i)=Pr[X_i=1]=\frac{r}{r^k}.
\end{equation*}

From the linear of the expectation \cite{alon},
\begin{eqnarray}\label{expectation}
E(X) &=& E(X_1+X_2+\cdots+X_q)\nonumber\\
     &=& \frac{r}{r^k}q \nonumber\\
     &=& \frac{n!}{(n-2k)!2^k k! r^{k-1}}.
\end{eqnarray}

By careful calculation, the following assertions (*) for
(\ref{expectation}) are true, which will be used later:
\begin{enumerate}
\item If $r$ is fixed, then for every $1\leq k\leq
\frac{n}{2}$, $E(X)\to\infty$.\\
\item If $k$ is fixed and $r\ll(\frac{n!}{(n-2k)! 2^k k!})^{\frac{1}{k-1}}$, then $E(X)\to\infty$; \\
\item If $k$ is fixed and $r\gg(\frac{n!}{(n-2k)! 2^k k!})^{\frac{1}{k-1}}$, then $E(X)\to 0 $; \\
\item If $k$ is fixed and $r=c^{(0)}(\frac{n!}{(n-2k)! 2^k
k!})^{\frac{1}{k-1}}$, where $c^{(0)}>0$ is a constant, then
$E(X)\to\frac{1}{(c^{(0)})^{k-1}}$.
\end{enumerate}

Though $k$ and $r$ can be both functions of $n$, if they are both
variables, the situation becomes very complicated. So we
illustrate monochromatic $k$-matching problem from three aspects:
$r$ is fixed, $k$ is fixed and $k=\lfloor\frac{n}{2}\rfloor$. The
last case is the perfect matching case or the nearly perfect
matching case. Since we focus on the asymptotic properties, we
will not distinguish $\lfloor\frac{n}{2}\rfloor$ from
$\frac{n}{2}$. That is, we always suppose $n$ is an even.

\subsection{$r$ is fixed}

Assertion (*) 1 says that $E(X)\to\infty$ for every $1\leq k\leq
\frac{n}{2}$ if $r$ is fixed. We certainly expect that $Pr[X>0]\to
1$ holds. In fact, it does.

\begin{theo}
If $r\geq 1$ is fixed, then almost every $\mathcal{K}_n^r$ has a
monochromatic $k$-matching for any $1\leq k\leq\frac{n}{2}$.
\end{theo}
\pf We have mentioned in Section 1 that the subgraph
$\mathcal{G}_i$ of $\mathcal{K}_n^r$ is actually the random graph
$\mathcal{G}_{n,p(n)}$, where $p(n)=\frac{1}{r}$. There is a
result saying that the threshold function for $\mathcal{G}_{n,p}$
has a perfect matching is $\frac{logn}{n}$ (\cite{jans}, page 85).
If $r$ is fixed, then $\frac{1}{r}\gg\frac{logn}{n}$, which
implies that almost every $\mathcal{G}_i$ has a perfect matching.
Then almost every $\mathcal{K}_n^r$ has a monochromatic
$k$-matching for every $1\leq k\leq\frac{n}{2}$. \qed

\subsection{$k$ is fixed}

In this case, we prove the following theorem.

\begin{theo}\label{monoma}
If $k$ is fixed ($k=1$ is a trivial case so suppose $k\geq 2$),
then
\begin{equation}
\lim_{n\to\infty} Pr[X>0]=\left\{
\begin{array}{ll}
0 & \mbox{ if }r\gg (\frac{n!}{(n-2k)!2^k k!})^{\frac{1}{k-1}},\\
1 & \mbox{ if }r\ll (\frac{n!}{(n-2k)!2^k k!})^{\frac{1}{k-1}}.
\end{array} \right.
\end{equation}

That is to say, $(\frac{n!}{(n-2k)!2^k k!})^{\frac{1}{k-1}}$ is
the threshold function for the property that $\mathcal{K}_n^r$ has
a monochromatic $k$-matching.
\end{theo}
\pf From {\it Markov's inequality} \cite{erdo}
$$Pr[X>0]\leq E(X)$$
and assertion (*) 3, we have
$$Pr[X>0]\to 0\mbox{ if }r\gg (\frac{n!}{(n-2k)!2^k k!})^{\frac{1}{k-1}}.$$

For the other half, we estimate $\frac{\Delta}{(E(X))^2}$, where
$\Delta=\sum_{i\sim j}Pr[A_i\cap A_j]$. $A_i(A_j)$ denotes the
event that the edges in $M_i(M_j)$ are monochromatic and $i\sim j$
means the ordered pair of $A_i$ and $A_j$ that are not independent
from each other.

Our goal is to prove that if $r\ll(\frac{n!}{(n-2k)!2^k
k!})^{\frac{1}{k-1}}$, then $\frac{\Delta}{(E(X))^2}\to 0.$
Because
\begin{eqnarray*}
\Delta &=& \sum_{i\sim j}Pr[A_i\cap A_j]\\
       &=& \sum_{s=1}^{k-1} \sum_{(i,j)_s}\frac{r}{r^{2k-s}}\\
       & & (\mbox{where }(i,j)_s\mbox{ means the ordered pair of }
       M_i \mbox{ and }M_j \mbox{ that have }s \mbox{ common
       edges})\\
       &\leq&\sum_{s=1}^{k-1}\frac{{n\choose 2}{{n-2}\choose 2}\cdots{{n-2(s-1)}\choose
       2}}{s!}(\frac{{{n-2s}\choose{2}}{{n-2s-2}\choose{2}}\cdots{{n-2k+2}\choose{2}}}{(k-s)!})^2\frac{1}{r^{2k-s-1}}\\
       &=&\frac{n!}{2^{2k}(n-2k)!(n-2k)!r^{2k-1}}\sum_{s=1}^{k-1}\frac{(n-2s)!2^sr^s}{s!(k-s)!(k-s)!},
\end{eqnarray*}
then we have
\begin{equation}\label{sec}
\frac{\Delta}{(E(X))^2}\leq
\frac{k!k!}{n!}\sum_{s=1}^{k-1}\frac{(n-2s)!2^sr^{s-1}}{s!(k-s)!(k-s)!}.
\end{equation}

If $r\ll(\frac{n!}{(n-2k)!2^k k!})^{\frac{1}{k-1}}\sim
(\frac{1}{2^kk!})^{\frac{1}{k-1}}n^\frac{2k}{k-1}$, then there are
3 possible cases: (i) $r\ll n^2$, (ii) $r=c^{(1)}n^2$, where
$c^{(1)}>0$ is a constant and (iii) $n^2\ll
r\ll(\frac{n!}{(n-2k)!2^k k!})^{\frac{1}{k-1}}$.

In case (i),
\begin{equation}\label{sum1}
\sum_{s=1}^{k-1}\frac{(n-2s)!2^sr^s}{s!(k-s)!(k-s)!}=(1+\circ(1))\frac{2(n-2)!}{(k-1)!(k-1)!}.
\end{equation}

Then submit (\ref{sum1}) to (\ref{sec}), we get
\begin{equation}\label{second1}
\frac{\Delta}{(E(X))^2} \leq 2(1+\circ(1))\frac{k^2}{n(n-1)}\to 0.
\end{equation}

In case (ii)
\begin{equation}\label{sum2}
\sum_{s=1}^{k-1}\frac{(n-2s)!2^sr^s}{s!(k-s)!(k-s)!}=c^{(2)}\frac{2(n-2)!}{(k-1)!(k-1)!},
\end{equation}
where $c^{(2)}$ is a sufficiently large constant.

Then submit (\ref{sum2}) to (\ref{sec}), we get
\begin{equation}\label{second2}
\frac{\Delta}{(E(X))^2} \leq 2c^{(2)}\frac{k^2}{n(n-1)}\to 0.
\end{equation}

In case (iii)
\begin{equation}\label{sum3}
\sum_{s=1}^{k-1}\frac{(n-2s)!2^sr^s}{s!(k-s)!(k-s)!}=(1+\circ(1))\frac{(n-2k+2)!2^{k-1}r^{k-2}}{(k-1)!}.
\end{equation}

Then submit (\ref{sum3}) to (\ref{sec}), we get
\begin{equation}\label{second3}
\frac{\Delta}{(E(X))^2} \leq
c^{(3)}\frac{n^{\frac{2k(k-2)}{k-1}}}{n^{2k-2}}\to 0,
\end{equation}
where $c^{(3)}$ is a sufficiently large constant.

Summarizing (\ref{second1}) (\ref{second2}) and (\ref{second3}),
we end the proof of $\frac{\Delta}{(E(X))^2}\to 0$ with the
condition $r\ll(\frac{n!}{(n-2k)!2^k k!})^{\frac{1}{k-1}}$.

A corollary of the {\it Chebyshev's inequality} \cite{erdo}
asserts that if $E(X)\to\infty$ and $\Delta=\circ((E(X))^2)$, then
almost surely $X>0$(\cite{alon},page 46). So from assertion (*) 2
and the above discuss, we obtain
$$Pr[X>0]\to 1\mbox{ if }r\ll (\frac{n!}{(n-2k)!2^k k!})^{\frac{1}{k-1}}.$$

From the definition of the threshold function (Definition
\ref{thre}), we can say that $(\frac{n!}{(n-2k)!2^k
k!})^{\frac{1}{k-1}}$ is the threshold function for the property
that $\mathcal{K}_n^r$ has a monochromatic $k$-matching. \qed

\subsection{$k=\frac{n}{2}$}
When $k=\frac{n}{2}$, a monochromatic $k$-matching is a
monochromatic perfect matching.

Replace $k$ with $\frac{n}{2}$ in (\ref{expectation}), we have
\begin{equation}\label{expectperf}
E(X)=\frac{n!}{(\frac{n}{2})!2^{\frac{n}{2}}r^{\frac{n}{2}-1}}.
\end{equation}
By calculation of (\ref{expectperf}), we get $E(X)\to 0$ if
$r\geq\frac{n }{c^{(4)}}$, where $c^{(4)}<e$ is a constant;
$E(X)\to \infty$ if $r\leq\frac{n }{e}$.

The following assertion is true as a direct corollary of Markov's
inequality and the threshold function for the property that
$\mathcal{G}_{n,p}$ having a perfect matching (\cite{jans}, page
85). Here we omit its proof.

\begin{theo}
If $r\geq\frac{n}{c^{(4)}}$, where $c^{(4)}<e$ is a constant, then
almost no $\mathcal{K}_n^r$ has a monochromatic perfect matching.
On the other hand, if $r\leq\frac{n}{logn+c^{(5)}(n)}$, where
$c^{(5)}(n)\to\infty$, then almost every $\mathcal{K}_n^r$ has a
monochromatic perfect matching.
\end{theo}

\section{Heterochromatic $k$-Matchings in $\mathcal{K}_n^r$}
Following the symbols in the previous section, let $B_i$ be the
event that the edges in $M_i$ are heterochromatic and $Y_i$ be the
indicator variable for the event $B_i$. That is,
\begin{equation}
Y_i= \left\{
\begin{array}{ll}
1 & \mbox{if } B_i \mbox{ happens},\\
0 & \mbox{otherwise}.
\end{array}
\right.
\end{equation}

Then for each $1\leq i\leq q$,
\begin{equation*}
Pr[Y_i=1]=\frac{{{r}\choose{k}}k!}{r^k}.
\end{equation*}

Then the random variable
$$Y=Y_1+Y_2+\cdots+Y_q$$
denotes the number of heterochromatic $k$-matchings in
$\mathcal{K}_n^r$.

From the linear of the expectation \cite{alon},
\begin{eqnarray}\label{expectation2}
E(Y) &=& E(Y_1+Y_2+\cdots+Y_q)\nonumber\\
     &=& \frac{r}{r^k}q \nonumber\\
     &=& \frac{n!}{(n-2k)!2^k k!}\frac{r!}{(r-k)!r^k}.
\end{eqnarray}

Since $r\geq k$ is a necessary condition in the heterochromatic
$k$-matching problem, we have the following assertion for $E(Y)$
by calculation of (\ref{expectation2}).

\begin{lemm}\label{sec3lem}
For every $1\leq k\leq n^{1-\epsilon}$ and $r\geq k$,
$E(Y)\to\infty$, where $0<\epsilon<1$ is a constant that can be
arbitrarily small.
\end{lemm}

The main result of this section is the following theorem:

\begin{theo}\label{3main}
If $1\leq k\leq n^{1-\epsilon}$ and $r\geq k$, where
$0<\epsilon<1$ is a constant that can be arbitrarily small, then
almost every $\mathcal{K}_n^r$ contains a heterochromatic
$k$-matching.
\end{theo}
\pf Similar to Theorem \ref{monoma}, for heterochromatic
$k$-matchings, the following estimate is for $\Delta'= \sum_{i\sim
j}Pr[B_i\cap B_j]$.

\begin{eqnarray*}
\Delta' &=& \sum_{i\sim j}Pr[B_i\cap B_j]\\
       &=& \sum_{s=1}^{k-1} \sum_{(i,j)_s}\frac{{r\choose k}k!{{r-s}\choose{k-s}}(k-s)!}{r^{2k-s}}\\
       & & (\mbox{where }(i,j)_s\mbox{ means the ordered pair of }
       M_i \mbox{ and }M_j \mbox{ that have }s \mbox{ common
       edges})\\
       &\leq& \sum_{s=1}^{k-1} \frac{{n\choose 2}{{n-2}\choose 2}\cdots{{n-2(s-1)}\choose 2}}{s!}
       (\frac{{{n-2s}\choose 2}{{n-2(s+1)}\choose 2}\cdots{{n-2(k-1)}\choose 2}}{(k-s)!})^2
       \frac{{r\choose k}k!{{r-s}\choose{k-s}}(k-s)!}{r^{2k-s}}\\
       &=&\frac{n!r!}{2^{2k}(n-2k)!(n-2k)!(r-k)!(r-k)!r^{2k}}\sum_{s=1}^{k-1}\frac{(n-2s)!(r-s)!2^sr^s}{s!(k-s)!(k-s)!}.
\end{eqnarray*}

Then
\begin{equation}\label{momen}
\frac{\Delta'}{(E(Y))^2}\leq\frac{k!k!}{n!r!}\sum_{s=1}^{k-1}\frac{(n-2s)!(r-s)!(2r)^s}{s!(k-s)!(k-s)!}.
\end{equation}

By careful calculation of (\ref{momen}), we get if $k\ll n$, then
\begin{eqnarray}\label{momen2}
\frac{\Delta'}{(E(Y))^2}&\leq&\frac{k!k!}{n!r!}(1+\circ(1))\frac{2r!(n-2)!}{(k-1)!(k-1)!}\nonumber\\
                        &=& (1+\circ(1))\frac{k^2}{n(n-1)}\to 0.
\end{eqnarray}

From (\ref{momen2}), Lemma \ref{sec3lem} and the assertion that if
$E(Y)\to\infty$ and $\Delta'=\circ((E(Y))^2)$, then almost surely
$Y>0$(\cite{alon}, page 46), we have
\begin{equation*}
Pr[Y>0]\to 1 ,
\end{equation*}
which finishes the proof. \qed

\begin{rema}
As a corollary of Theorem \ref{3main}, if one of $k$ and $r(\geq
k)$ is fixed, then almost every $\mathcal{K}_n^r$ has a
heterochromatic $k$-matching. The only left case that we can not
deal with is that $k=c^{(6)}n$, where $0<c^{(6)}\leq\frac{1}{2}$
is a constant.
\end{rema}
\section{Results on Other Subgraphs}

Completely similar to Section 2 and Section 3, we can study
monochromatic $k$-clique, $k$-tree and heterochromatic $k$-clique,
$k$-tree in $\mathcal{K}_n^r$. We list our results here.

\begin{theo}
If $r$ is fixed, then
\begin{equation*}
\lim_{n\to\infty}Pr[\mathcal{K}_n^r \mbox{ contains a
monochromatic $k$-clique}]=\left\{
\begin{array}{ll}
0 & \mbox{ if }k\geq 2log_rn,\\
1 & \mbox{ if }k\leq \frac{log_rn}{1.704\times 10^9}.
\end{array} \right.
\end{equation*}
\end{theo}

\begin{theo}
If $k$ is fixed, then
\begin{equation*}
\lim_{n\to\infty}Pr[\mathcal{K}_n^r \mbox{ contains a
monochromatic $k$-clique}]=\left\{
\begin{array}{ll}
0 & \mbox{ if }r\gg n^{\frac{k}{{k\choose 2}-1}},\\
1 & \mbox{ if }r\leq (\frac{1}{2k!})^{\frac{1}{{k\choose
2}-1}}n^{\frac{k}{{k\choose 2}-1}}.
\end{array} \right.
\end{equation*}
That is to say, $n^{\frac{k}{{k\choose 2}-1}}$ is the threshold
function for the property that $\mathcal{K}_n^r$ has a
monochromatic $k$-clique.
\end{theo}

\begin{theo}
If $r\geq n^{4+\epsilon}$, where $\epsilon>0$ is a constant that
can be arbitrarily small, then for every $k\leq n$, there almost
surely exists a heterochromatic $k$-clique in $\mathcal{K}_n^r$.
\end{theo}

\begin{theo}
If $k$ is fixed, then
\begin{equation*}
\lim_{n\to\infty}Pr[\mathcal{K}_n^r \mbox{ contains a
monochromatic $k$-tree}]=\left\{
\begin{array}{ll}
0 & \mbox{ if }r\gg k{n\choose k}^{\frac{1}{k-2}},\\
1 & \mbox{ if }r\leq \frac{k}{n}{n\choose k}^{\frac{1}{k-2}}.
\end{array} \right.
\end{equation*}
\end{theo}

\begin{theo}
If $r\geq c^{(7)}n$, where $c^{(7)}>1$ is a constant, then almost
no $\mathcal{K}_n^r$ contains a monochromatic spanning tree.
\end{theo}

\begin{theo}
If $r$ is fixed, then almost every $\mathcal{K}_n^r$ contains a
monochromatic $k$-tree for any $2\leq k\leq n$.
\end{theo}

\begin{theo}
If $2\leq k\leq logn$ and $r\geq k-1$, then almost every
$\mathcal{K}_n^r$ contains a heterochromatic $k$-tree.
\end{theo}


\begin{thebibliography}{99}

\bibitem{alon}
N. Alon and J. Spencer, {\it The Probabilistic Method, 2nd ed.,}
John Wiley \& Sons, Inc. 2000.

\bibitem{bond}
J. A. Bondy and U. S. R. Murty, {\it Graph Theory with
Applications,} The Macmillan Press LTD., 1976.

\bibitem{boll}
B. Bollob\'{a}s, {\it Random Graphs, 2nd ed.,} Cambridge
University Press, 2001.

\bibitem{erdo}
P. Erd\"{o}s and J. Spencer, {\it Probabilistic Methods in
Combinatorics,} Academic Press, 1974.

\bibitem{spen}
J. Spencer, {\it The Strange Logic of Random Graphs,} Springer,
2001.

\bibitem{jans}
S. Janson, T. \L uczak and A. Rucinski, {\it Random Graphs}, John
Wiley \& Sons, Inc., 2000.

\bibitem{erdo1}
P. Erd\"{o}s and A. R\'{e}nyi, On random graphs I, {\it Publ.
Math. Debrecen} {\bf 6}, 290-297.

\bibitem{erdo2}
P. Erd\"{o}s and A. R\'{e}nyi, On the evolution of random graphs,
{\it Publ. Math. Inst. Hungar. Acad. Sci.} {\bf 5}, 17-61.

\bibitem{erdo3}
P. Erd\"{o}s and A. R\'{e}nyi, On the evolution of random graphs,
{\it Bull. Inst. Int. Statist. Tokyo} {\bf 38}, 343-347.

\end{thebibliography}
\end{document}